\documentclass{amsart}
\usepackage{amsmath,amssymb}

\newtheorem{conjecture}{Conjecture}

\usepackage{epsfig}
\usepackage{url}
\usepackage{amssymb}
\usepackage{amsmath}
\usepackage{amsfonts}

\def\Dbar{\leavevmode\lower.6ex\hbox to 0pt{\hskip-.23ex \accent"16\hss}D}

\def\bZ{{\mbox{\bf Z}}}

\begin{document}

\title{A few new orders for D-optimal matrices}
\author {Dragomir {\v{Z}}. {\Dbar}okovi{\'c}}
\address{University of Waterloo, 
Department of Pure Mathematics and Institute for Quantum Computing,
Waterloo, Ontario, N2L 3G1, Canada}
\email{djokovic@uwaterloo.ca}
\date{}

\begin{abstract}
The first examples of D-optimal matrices of orders $222$, $234$, $258$ and $278$ are constructed.
\end{abstract}

\maketitle

\section{Introduction}
\label{uvod}

A square $\{\pm 1\}$-matrix $M$ of order $m$ is {\em D-optimal} if its determinant is
maximal among all such matrices. It is customary to denote the maximal value of these determinants 
by $\alpha_m$. 
D-optimal matrices are also known as {\em D-optimal designs}.
If $m$ is divisible by 4 or $m<3$ then the D-optimal matrices of order $m$
are just the Hadamard matrices of order $m$ with positive determinant. In this note we assume from now on 
that $m\equiv 2 \pmod {4}$ and we set $m=2v$ where $v$ is odd.
Then it is well-known that $\alpha_m\le 2^v (m-1)(v-1)^{v-1}$ and that this inequality is strict 
if $m-1$ is not a sum of two squares. For this and other well-known facts see e.g. \cite{Brent:2013,KO-Dopt:2007}.
When $m-1$ is a sum of two squares we expect that the equality will hold (see Conjecture 1 below).

We use the well known $2 \times 2$ array
\begin{equation} 
M=\left[ \begin{array}{cc}
A    & B \\
-B^T & A^T 
\end{array} \right] \label{mat-M}
\end{equation}
in order to construct our new D-optimal matrices. Here, $A$ and $B$ are $\{\pm 1\}$-matrices of order $v$
and the superscript T denotes transposition. It is well-known that $M$ is D-optimal
if $A$ and $B$ commute and satisfy the equations
\begin{equation} \label{ABjedn}
AA^T + BB^T=A^T A+B^T B= (m-2)I_v+2J_v,
\end{equation}
where $I_v$ is the identity matrix and $J_v$ the matrix with all entries equal to $+1$.
The matrices $A$ and $B$ are usually constructed from difference families (DF) consisting
of two blocks $X$ and $Y$ in some abelian group $G$ of order $v$.
(We do not require that $X$ and $Y$ have the same size.)
A DF is {\em cyclic} if $G$ is a cyclic group.

Let $(v; r=|X|,s=|Y|; \lambda)$ be the parameter set (PS) of such a DF, and define an additional
parameter $n$ by $n=r+s-\lambda$. We say that such a DF and its PS are {\em D-optimal}
if $2v-1$ is a sum of two squares and $v=2n+1$.

If $(X,Y)$ is a D-optimal DF then we can interchange $X$ and $Y$ or replace $X$ by its complement in $G$ and we obtain
a new D-optimal DF having the same parameter $n$. Therefore we may (and we do) assume that $(v-1)/2 \ge r\ge s$.
We say that the PSs satisfying these inequalities are {\em normalized}.
The normalized D-optimal PSs are in one-to-one correspondence with integer pairs $(x,y)$ such that $x\ge y\ge 0$.
(For convenience we do not omit the three trivial cases with $1 \ge x \ge y \ge 0$.)
Given such a pair, the corresponding normalized D-optimal PS is given by the formulas 
\begin{eqnarray}
	v &=& 1+x+x^2+y+y^2, \label{par-v} \\
	r &=& x(x+1)/2+y(y-1)/2, \label{par-r} \\
	s &=& x(x-1)/2+y(y+1)/2, \label{par-s} \\
	\lambda &=& x(x-1)/2+y(y-1)/2, \label{par-l} 
\end{eqnarray}
see \cite[Proposition 1]{DK:2015}. Conversely, if $(v;r,s;\lambda)$ is a normalized D-optimal PS then the
corresponding pair $(x,y)$ is given by $x=r-\lambda$, $y=s-\lambda$.

Each D-optimal DF $(X,Y)$ gives rise to a D-optimal matrix.
Briefly, the blocks $X$ and $Y$ provide commuting $\{\pm 1\}$-matrices $A$ and $B$
of order $v$ satisfying the equation \eqref{ABjedn}. The corresponding D-optimal matrix $M$ is obtained 
by plugging $A$  and $B$ into the array \eqref{mat-M}.
When this D-optimal DF is cyclic, say $G=\bZ_v$ (integers mod $v$), $A$ and $B$ are circulants and one obtains
the first row, say $a_0,a_1,\ldots,a_{v-1}$, of $A$ by setting $a_i=-1$ if $i\in X$ and $a_i=+1$ otherwise.

The normalized D-optimal PSs with $v<100$ are listed in \cite{DK:2015} and for each of them
(with two exceptions) one cyclic D-optimal DF is given.
The two exceptional cases are $(85;39,34;31)$ and $(99; 43,42; 36)$ and it is still unknown
wheather the corresponding DFs exist.
In the next section we list all normalized D-optimal PSs in the range $100<v<200$ and give the references for the known
D-optimal DFs for them (including our new results).

Apparently the following conjecture has not been proposed so far.
\begin{conjecture}
	For each D-optimal parameter set $(v;r,s;\lambda)$ there exists a difference family
	(over some abelian group of order $v$) having these parameters.
\end{conjecture}
No counterexamples are known even for cyclic DFs. The smallest undecided case is $(85; 39,34; 31)$.

Among all odd integers $v<200$, for which $2v-1$ is a sum of two squares, no D-optimal matrices
of order $2v$ are known only in the following 17 cases:
\begin{equation*} 
v = 99,115,123,135,141,147,153,159,163,167,169,175,177,185,187,195,199.
\end{equation*}

Our main result is the construction of the first examples of D-optimal matrices of size $2v$ for 
$v=111,117,129,139$. They are listed in section \ref{glavni}. In some cases we give several 
nonequivalent solutions. In the case $v=111$ there are two normalized D-optimal PSs and 
we provide five solutions to each. In the other three cases there is only one normalized PS
for each of them. In the next section we list the normalized D-optimal PSs in the range $100 < v < 200$
and provide references if the solutions are known.

There are two borderline cases: the case where $y=0$ and the case $x=y$.
Note that $y=0$ is equivalent to $\lambda=s$, and $x=y$ is equivalent to $r=s$.
We consider these borderline cases in sections \ref{prv-sl} and \ref{drg-sl}.

\section{Normalized D-optimal parameter sets with $100<v<200$}
\label{100-200}

There are 40 normalized D-optimal parameter sets with $100<v<200$. They are listed in Table 1. 
If a cyclic D-optimal design is known, we give at least one reference where it occurs. 
The question mark means that, for the $v$ given in the parameter set, we are not aware 
that any D-optimal matrix of order $2v$ is known.
The asterisk means that the required cyclic D-optimal design is constructed in this note. 

The group of units $\bZ^*_v$ acts on $\bZ_v$ by multiplication.
The same is true for any subgroup $H$ of $\bZ_v^*$.
We choose a suitable $H$ and construct a D-optimal DF whose blocks $X$ and $Y$
can be represented as a union of $H$-orbits.
We specify $X$ and $Y$ by writing them as $HR$ where $R$ is a set of representatives of
the $H$-orbits contained in the corresponding block.

\vspace{3 mm}
\begin{tabular}{cl|cl}
	\multicolumn{4}{c}{Table 1: Normalized D-optimal PSs for $100<v<200$} \\
	\hline
	(103; 46,43; 38) & \cite{DK:JCD:2012} &            (153; 70,66; 60) & ? \\
	(103; 48,42; 39) & \cite{DK:JCD:2012} &            (153; 72,65; 61) & ? \\
	(111; 51,46; 42) & * &                             (157; 78,66; 66) & \cite{Gysin:PhD:1997,GS:JCMCC:1998} \\
	(111; 55,45; 45) & * &                             (159; 78,67; 66) & ? \\
	(113; 49,49; 42) & \cite{Gysin:PhD:1997,Gysin:AJC:1997} & (163; 73,72; 64) & ? \\
	(113; 55,46; 45) & ? &                             (163; 76,70; 65) & ? \\
	(115; 51,49; 43) & ? &                             (163; 79,69; 67) & ? \\
	(117; 56,48; 46) & * &                             (167; 76,73; 66) & ? \\
	(121; 55,51; 46) & \cite{DK:JCD:2012} &            (169; 81,72; 69) & ? \\
	(123; 58,51; 48) & ? &                             (175; 81,76; 70) & ? \\
	(129; 57,56; 49) & * &                             (177; 84,76; 72) & ? \\
	(131; 61,55; 51) & \cite{DK:JCD:2012} &            (181; 81,81; 72) & \cite{Gysin:PhD:1997,GS:JCMCC:1998} \\
	(133; 60,57; 51) & ? &                             (183; 83,81; 73) & ? \\
	(133; 66,55; 55) & \cite{KKS:1991} &               (183; 91,78; 78) & \cite{KKS:1991} \\
	(135; 66,56; 55) & ? &                             (185; 91,79; 78) & ? \\
	(139; 67,58; 56) & * &                             (187; 88,81; 76) & ? \\
	(141; 65,60; 55) & ? &                             (189; 87,83; 76) & ? \\
	(145; 64,64; 56) & \cite{DZD:AJC:1997,Gysin:PhD:1997,GS:JCMCC:1998,KSWX:1997} & (189; 92,81; 79) & ? \\
	(145; 69,61; 58) & ? &                             (195; 94,84; 81) & ? \\
	(147; 66,64; 57) & ? &                             (199; 93,87; 81) & ? \\ 
	\hline
\end{tabular} 
\vspace{3 mm}

There are two known infinite series of D-optimal matrices.

The first one, constructed by C. Koukouvinos, S. Kounias and J. Seberry \cite{KKS:1991},
consists of cyclic D-optimal designs with parameter sets
$$(1+q+q^2; q(q+1)/2,q(q-1)/2; q(q-1)/2)$$
where $q$ is any prime power.

The second one consists of D-optimal matrices of order $2(1+2q+2q^2)$ where $q$ is any odd prime power.
L. E. Brouwer constructed an infinte family of symmetric designs (also known as symmetric balanced incomplete 
block designs or SBIBD) with parameters 
$$ v=1+2(q+q^2+\cdots+q^h),\quad k=q^h,\quad \lambda=(q^h-q^{h-1})/2 $$
where $q$ is an odd prime power and $h\ge1$ an integer.
The subfamily obtained by setting $h=2$ was used by A. L. Whiteman \cite{Whiteman:1990}
to construct an infinte series of D-optimal matrices of order mentioned above.
No DFs are used in his construction. 
More precisely, if $S$ is the $\{\pm1\}$-incidence matrix of a SBIBD with parameters
$v=1+2q+2q^2$, $k=q^2$, $\lambda=q(q-1)/2$ then the matrix
$$\left[ \begin{array}{cc}
S  & S \\
-S & S 
\end{array} \right]$$
is a D-optimal matrix of order $2v$.

We shall refer to these two infinite series as the {\em KKS series} and the {\em Whiteman series}, respectively.

\section{D-optimal matrices of orders 222, 234, 258 and 278}
\label{glavni}

No D-optimal matrices of these orders were known previously. We construct them 
from cyclic D-optimal DFs $(X,Y)$ for $v=111, 117, 129, 139.$
These DFs are presented below in the following way. First we give the 
parameter set $(v; r,s; \lambda)$ and the subgroup $H$ of $\bZ^*_v$ that we use. 
The blocks $X$ and $Y$ should be computed by using the formulas
$X=\cup_{i\in I} Hi$ and $Y=\cup_{j\in J} Hj$,
where $Hi$ is the orbit of $H$ containing the integer $i\in\bZ_v$.
In particular we have $H0=\{0\}$.  The sets $I$ and $J$ are representatives of the
$H$-orbits contained in $X$ and $Y$, respectively.
The numeric subscripts attached to $I$ and $J$ are used to label 
different D-optimal DFs for the same parameter set.
All DFs listed for the same PS are pairwise nonequivalent. 

For $v=111$ there are two D-optimal PSs and for each of them we computed five DFs:

\begin{eqnarray*}
&& \qquad \qquad (111; 51,46; 42), \quad H=\{1,10,100\} \\
I_1 &=& \{3,8,13,14,15,16,21,26,31,32,43,51,52,53,55,63,64\} \\
J_1 &=& \{0,2,3,5,7,13,14,15,26,31,43,52,53,54,63,64\} \\
\\
I_2 &=& \{1,4,8,13,15,21,22,26,27,41,42,43,44,52,53,63,64\} \\
J_2 &=& \{0,1,2,4,11,13,14,17,21,25,33,42,52,53,62,63\} \\
\\
I_3 &=& \{1,3,5,7,9,13,14,15,17,27,31,42,44,51,52,62,63\} \\
J_3 &=& \{0,1,2,4,11,17,21,25,42,43,51,53,54,55,62,64\} \\
\\
I_4 &=& \{1,2,6,8,13,14,15,16,21,27,31,41,43,52,53,55,63\} \\
J_4 &=& \{0,1,2,4,5,11,16,17,21,22,32,33,43,44,54,63\} \\
\\
I_5 &=& \{1,2,4,8,9,15,16,26,31,33,42,44,51,53,54,55,63\} \\
J_5 &=& \{0,2,3,4,8,13,16,21,41,44,51,52,55,62,63,64\}
\end{eqnarray*}

\begin{eqnarray*}
&& \qquad \qquad (111; 55,45; 45), \quad H=\{1,10,100\} \\
I_1 &=& \{0,2,6,13,14,16,17,18,22,26,33,43,51,52,53,54,55,63,64\} \\
J_1 &=& \{2,3,4,6,13,14,17,21,25,27,41,54,55,63,64\} \\
\\
I_2 &=& \{0,3,4,6,7,8,9,13,16,25,32,41,42,44,52,54,62,63,64\} \\
J_2 &=& \{2,3,5,7,13,17,21,25,41,44,51,52,55,62,63\} \\
\\
I_3 &=& \{0,5,6,9,13,14,17,21,25,31,32,33,41,42,43,53,62,63,64\} \\
J_3 &=& \{2,5,6,8,13,16,18,31,33,43,51,54,55,62,63\} \\
\\
I_4 &=& \{0,3,4,5,7,14,15,18,22,26,27,32,42,43,52,55,62,63,64\} \\
J_4 &=& \{2,5,9,15,18,22,25,32,33,41,42,43,44,52,53\} \\
\\
I_5 &=& \{0,3,4,5,7,9,15,16,22,25,26,32,42,44,51,53,54,55,64\} \\
J_5 &=& \{2,4,14,15,17,18,21,25,41,43,51,52,54,63,64\}
\end{eqnarray*}

\newpage
For $v=117$ we computed eight D-optimal DFs: 

\begin{eqnarray*}
&& \qquad \qquad (117; 56,48; 46), \quad H=\{1,16,22\} \\
I_1 &=& \{2,7,8,9,10,13,18,19,20,39,41,47,51,56,58,63,73,78,79,95\} \\
J_1 &=& \{0,1,2,5,6,8,12,19,25,34,39,41,47,56,57,63,78,79\} \\
\\
I_2 &=& \{3,5,6,7,18,19,21,24,25,26,28,29,39,42,47,56,57,58,78,95\} \\
J_2 &=& \{0,1,3,7,14,24,25,29,34,39,41,47,51,56,57,63,78,79\} \\
\\
I_3 &=& \{1,3,7,13,17,19,20,24,26,28,29,39,41,47,51,57,58,78,79,95\} \\
J_3 &=& \{0,4,7,9,12,13,17,18,20,25,28,29,34,35,39,51,56,78\} \\
\\
I_4 &=& \{1,2,3,4,5,6,10,12,13,14,19,20,24,29,35,39,47,73,78,79\} \\
J_4 &=& \{0,4,6,13,14,21,26,28,29,36,39,40,56,57,58,73,78,95\} \\
\\
I_5 &=& \{2,3,5,7,10,13,21,24,25,26,28,35,39,42,51,56,57,58,63,78\} \\
J_5 &=& \{0,3,5,13,14,17,19,20,25,34,36,39,47,51,56,63,73,78\} \\
\\
I_6 &=& \{3,5,9,10,13,18,20,21,25,34,35,36,39,41,51,56,58,63,78,95\} \\
J_6 &=& \{0,8,9,10,12,13,14,19,26,28,35,39,51,56,57,58,63,78\} \\
\\
I_7 &=& \{3,4,5,6,7,12,13,17,19,24,26,36,39,41,42,47,57,73,78,79\} \\
J_7 &=& \{0,3,7,13,20,25,28,34,35,36,39,47,51,56,57,73,78,95\} \\
\\
I_8 &=& \{1,2,6,7,9,12,17,18,20,25,28,36,39,41,51,56,58,63,78,79\} \\
J_8 &=& \{0,1,2,5,9,12,13,18,28,29,39,40,41,47,51,56,73,78\}
\end{eqnarray*}

For $v=129$ we just have one D-optimal DF:

\begin{eqnarray*}
&& \qquad \qquad (129; 57,56; 49), \quad H=\{1,49,79\} \\
I_1 &=& \{1,6,7,17,20,21,22,26,31,35,39,42,44,57,60,62,63,68,73\} \\
J_1 &=& \{0,1,3,11,12,19,20,22,31,35,39,44,50,63,65,68,70,73,78,86\}
\end{eqnarray*}

For $v=139$ we found only two D-optimal DFs:

\begin{eqnarray*}             
&& \qquad \qquad (139; 67,58; 56), \quad H=\{1,42,96\} \\
I_1 &=& \{0,3,5,8,9,15,21,22,23,24,25,34,39,41,46,49,62,66,69,72,75,82,85\} \\
J_1 &=& \{0,4,6,8,9,11,12,18,23,26,33,34,36,39,41,49,59,65,66,78\} \\
\\
I_2 &=& \{0,2,4,11,13,14,17,21,22,25,26,28,34,41,43,56,59,62,65,66,72,82,85\} \\
J_2 &=& \{0,1,3,11,12,13,15,17,21,23,25,26,31,33,36,41,65,68,69,75\}
\end{eqnarray*}

\section{The borderline case $\lambda=s$} \label{y=0}
\label{prv-sl}

In this case $y=0$ and the equations \eqref{par-v}-\eqref{par-l} imply that 
these borderline cases have parameters
\begin{equation} \label{PS-0}
(v=1+x+x^2; r=x(x+1)/2,s=x(x-1)/2; \lambda=s).
\end{equation}
They are listed in Table 2 for $x\le 15$ ($x=0$ is omitted).
The whole KKS series is contained in this borderline case.
The symbol "cyc" means that there exists a cyclic D-optimal DF $(X,Y)$ with
the PS given in the second column. We recall that all DFs in the KKS series
are also cyclic.
The question mark in Table 2 indicates that no D-optimal matrix of order 422 is known,
and the asterisk means that for $x=10$ the D-optimal DFs with parameters \eqref{PS-0}
are constructed in this note (see section \ref{glavni}).

\vspace{3 mm}
\begin{tabular}{|r|l|l|l|}
	\multicolumn{4}{c}{Table 2: Normalized D-optimal PSs with $\lambda=s$} \\
	\hline
	$x$ & parameter set        & status & references \\
	\hline
	1 & $(3; 1,0; 0)$          & KKS    & \cite{KSWX:1997} \\
	2 & $(7; 3,1; 1)$          & KKS    & \cite{KSWX:1997} \\
	3 & $(13; 6,3; 3)$         & KKS    & \cite{KSWX:1997} \\
	4 & $(21; 10,6; 6)$        & KKS    & \cite{KSWX:1997} \\
	5 & $(31; 15,10; 10)$      & KKS    & \cite{KSWX:1997} \\
	6 & $(43; 21,15; 15)$      & cyc    & \cite{Cohn:1989} \\
	7 & $(57; 28,21; 21)$      & KKS    & \cite{KSWX:1997} \\
	8 & $(73; 36,28; 28)$      & KKS    & \cite{KSWX:1997} \\
	9 & $(91; 45,36; 36)$      & KKS    & \cite{KSWX:1997} \\
	10 & $(111; 55,45; 45)$    & cyc    & ${*}$ \\
	11 & $(133; 66,55; 55)$    & KKS    & \cite{KSWX:1997} \\
	12 & $(157; 78,66; 66)$    & cyc    & \cite{Gysin:PhD:1997} \\
	13 & $(183; 91,78; 78)$    & KKS    & \cite{KSWX:1997} \\
	14 & $(211; 105,91; 91)$   & ?      & \\ 
	15 & $(241; 120,105; 105)$ & cyc    & \cite{DK:JCD:2012} \\
	\hline
\end{tabular} 
\vspace{3 mm}

Apart from the two infinite series of D-optimal matrices (the KKS and the Whiteman 
series) there are only finitely many known D-optimal matrices. The largest 
order among them is 482. It is achieved by the D-optimal matrix provided 
by the last case $(x=15)$ of Table 2.

\section{The borderline case $r=s$} \label{x=y}
\label{drg-sl}

In this case $y=x$ and the equations \eqref{par-v}-\eqref{par-l} imply that the PS is 
\begin{equation} \label{PS=}
(v=1+2x+2x^2; r=x^2,s=x^2; \lambda=x^2-x).
\end{equation}
These PSs are listed in Table 3 for $x\le15$ ($x=0$ is omitted).
Since $r=s$ the sets $X$ and $Y$ have the same cardinality.
If $Y=\mu X$ for some $\mu\in\bZ_v^*$ we say that $\mu$ is a {\em multiplier} of this DF.
Note that if $\mu$ is a multiplier of $(X,Y)$ then each element of the coset $\mu H$
is also a multiplier of $(X,Y)$.  Note also that if $\mu$ is a multiplier of a
DF $(X,Y)$ then $\mu^{-1}$ is a multiplier of the DF $(Y,X)$.
When we list two values for $\mu$, this means that for the given PS there exist
two nonequivalent DFs, one for each value of $\mu$.
For instance, for $x=2$ we list two values 1 and 5 for $\mu$. For $\mu=5$ we can take 
$X=\{0,1,3,9\}$ and $Y=\{0,2,5,6\}$ (see \cite{KSWX:1997}) while for $\mu=1$ 
we can take $X=Y=\{0,1,4,6\}$ (a difference set with parameters $(13,4,1)$).
In connection with this example we raise the following question. Does there exist 
a difference set (over some finite abelian group) with parameters
$v=1+2x+2x^2$, $k=x^2$, $\lambda=x(x-1)/2$ where $x>2$ is an integer?

The letter "W" in the third column means that the Whiteman series provides an
example of a D-optimal matrix of order $2v$.
The question mark indicates that no D-optimal matrix of order $2v$ is known.
The asterisk in the last column indicates that an example is provided below in this section.

\vspace{3 mm}
\begin{tabular}{|r|l|l|r|l|}
	\multicolumn{5}{c}{Table 3: Normalized D-optimal PSs with $r=s$} \\
	\hline
	$x$ & parameter set        & status & $\mu$  & references \\
	\hline
	1 & $(5; 1,1; 0)$          & W,cyc  & 3      & \\
	2 & $(13; 4,4; 2)$         & W,cyc  & 1,5    & \cite{KSWX:1997} \\
	3 & $(25; 9,9; 6)$         & W,cyc  & 7      & \cite{KSWX:1997} \\
	4 & $(41; 16,16; 12)$      & cyc    & 9      & \cite{Djokovic:RadMat:1991,KSWX:1997} \\
	5 & $(61; 25,25; 20)$      & W,cyc  & 11     & \cite{Djokovic:RadMat:1991,Gysin:PhD:1997,KSWX:1997} \\ 
	6 & $(85; 36,36; 30)$      & cyc    & 3      & ${*}$ \\
	7 & $(113; 49,49; 42)$     & W,cyc  & 2,15   & ${*}$, \cite{Gysin:PhD:1997,Gysin:AJC:1997,KSWX:1997} \\ 
	8 & $(145; 64,64; 56)$     & cyc    & 11,14  & ${*}$, \cite{Djokovic:RadMat:1991} \\
	9 & $(181; 81,81; 72)$     & W,cyc  & 19     & \cite{Gysin:PhD:1997,KSWX:1997} \\ 
	10 & $(221; 100,100; 90)$  & ?      &        & \\
	11 & $(265; 121,121; 110)$ & W      &        & \\
	12 & $(313; 144,144; 132)$ & ?      &        & \\
	13 & $(365; 169,169; 156)$ & W      &        & \\
	14 & $(421; 196,196; 182)$ & ?      &        & \\ 
	15 & $(481; 225,225; 210)$ & ?      &        & \\ 
	\hline
\end{tabular} 
\vspace{3 mm}

It is conjectured in \cite{KSWX:1997} that whenever $x$ is a prime power
then there exists a cyclic DF with the PS \eqref{PS=} and the multiplier
$\mu=2x+1$. The prime powers $x=2,3,4,5,7,9$ were given as examples.
The case $x=8$ was skipped. We found solutions for $x=6$ (not a prime power) and
$x=8$ but not with multiplier $2x+1$. 
We suggest that the conjecture we just cited may be false for $x=8$.

For $x=6$ we use the cyclic subgroup $H=\{1,9,16,19,21,49,59,81\}$ 
of order 8 of $\bZ_{85}^*$, 
The block $X$ consists of six $H$-orbits with representatives 
3,12,15,17,24,34. The second block is $Y=3X$. We have $|X|=|Y|=36$ and 
$|X\cap Y|=4$. Evidently, $\mu=3$ is a multiplier of $(X,Y)$.

For $x=7$ we use the subgroup $H=\{1,16,28,30,49,106,109\}$ of $\bZ_{113}^*$.
Let $X=HR$ where $R=\{1,2,3,7,17,27,40\} \subseteq \bZ_{113}^*$. Then $|X|=49$
and $(X,2X)$ is a D-optimal DF with multiplier 2. This example is not equivalent to 
the one in \cite{KSWX:1997} which has multiplier 15.

For $x=8$ we use the subgroup $H=\{1,16,36,81,111,136,141\}$ of $\bZ_{145}^*$.
Let $X=HR$ where $R=\{1,3,4,7,8,10,12,26,29,46\} \subseteq \bZ_{145}$.
The point $29\in\bZ_{145}$ is a fixed point of $H$ and all the other nine $H$-orbits
in $X$ have size 7. Hence $|X|=64$. Further, $(X,11X)$ and $(X,14X)$ are two 
nonequivalent D-optimal DFs with multipliers $11$ and $14$, respectively.
We also note that $|X\cap 11X|=15$ while $|X\cap 14X|=21$.

\section{Acknowledgements}
This research was enabled in part by support provided by the Shared Hierarchical 
Academic Research Computing Network (SHARCNET: www.sharcnet.ca) and
the Digital Research Alliance of Canada (alliancecan.ca).


\begin{thebibliography}{99}

\bibitem{Brent:2013}
R. P. Brent, Finding D-optimal designs by randomised decomposition and switching,
Australasian J. Combin. 55 (2013), 15--30. 

\bibitem{AEB:1983}
A. E. Brouwer, An infinite series of symmetric designs,
Mathematisch Centrum, Amsterdam, Report ZW 202/83 (1983).

\bibitem{Cohn:1989}
J. H. E. Cohn, 
On determinants with elements $\pm 1$, II.
Bull. London Math. Soc. 21 (1989), 36--42.

\bibitem{Djokovic:RadMat:1991}
D. {\v{Z}}. {\Dbar}okovi{\'c},
On maximal $\{+1,-1\}$-matrices of order $2n$, $n$ odd,
Radovi Matemati{\v c}ki 7 (1991), 371--378.

\bibitem{DZD:AJC:1997}
D. {\v{Z}}. {\Dbar}okovi{\'c},
Some new D-optimal designs. 
Austarlasian J. Combin. 15 (1997), 221–-231.

\bibitem{DK:JCD:2012}
D. {\v{Z}}. {\Dbar}okovi{\'c}, I. S. Kotsireas,
New results on D-optimal matrices.
J. Combin. Designs, 20 (2013), 278--289.

\bibitem{DK:2015}
D. {\v{Z}}. {\Dbar}okovi{\'c}, I. S. Kotsireas, D-optimal matrices of orders 
118, 138, 150, 154 and 174. In: Colbourne, C.J. (ed), 
Algebraic Design Theory and Hadamard Matrices,
Springer International Publishing Switzerland 2015, p. 71--82.
doi:10.1007/978-3-319-17729-8

\bibitem{DK:CompMethods:2019}
D. {\v{Z}}. {\Dbar}okovi{\'c}, I. S. Kotsireas,
Computational methods for difference families in finite abelian groups,
Spec. Matrices 2019; 7:127--141. 

\bibitem{Gysin:PhD:1997}
Marc Michel Gysin, Combinatorial designs, sequences and
cryptography, PhD Thesis, University of Wollonglong, 1997.

\bibitem{Gysin:AJC:1997}
M. Gysin,
New D-optimal designs via cyclotomy and generalized cyclotomy,
Australasian J. Combin. 15 (1997), 247--255.

\bibitem{GS:JCMCC:1998}
M. Gysin, J. Seberry,
An experimental search and new combinatorial designs via a generalization of cyclotomy,
JCMCC 27 (1998), p. 143--160.

\bibitem{KO-Dopt:2007}
H. Kharaghani and W. Orrick,
D-optimal matrices, in Handbook of Combinatorial Designs, 2nd ed.
C. J. Colbourn, J. H. Dinitz (eds)  pp. 296--298.
Discrete Mathematics and its Applications (Boca Raton).
Chapman \& Hall/CRC, Boca Raton, FL, 2007.

\bibitem{KKS:1991}
C. Koukouvinos, S. Kounias and J. Seberry, 
Supplementary diference sets and optimal designs, 
Discrete Math. 88 (1991), 49--58.

\bibitem{KSWX:1997}
C. Koukouvinos, J. Seberry, A. L. Whiteman, Ming-yuan Xia,
Optimal designs, supplementary difference sets and multipliers,
Journal of Statistical Planning and Inference 62 (1997) 81--90.

\bibitem{Whiteman:1990}
A. L. Whiteman, A family of D-optimal designs, Ars Combinatoria 30 (1990), 23--26.

\end{thebibliography}
\end{document}